\documentclass[11pt]{amsart}
\usepackage{amsmath,amssymb}
\usepackage{amscd,eucal,amsthm}
\newtheorem{Theorem}{Theorem}
\newtheorem{Corollary}{Corollary}
\newtheorem{Lemma}{Lemma}
\newtheorem{Fact}{Fact}
\newtheorem{Proposition}{Proposition}

\theoremstyle{remark}
\newtheorem*{Remark}{Remark}

\renewcommand{\to}[1][]{\xrightarrow{#1}}

\renewcommand{\gg}{{\mathfrak{g}}}

\newcommand{\nn}{{\mathfrak{n}}}
\newcommand{\z}{{\mathfrak{z}}}

\newcommand{\hh}{{\mathfrak{h}}}
\renewcommand{\l}{{\lambda}}

\newcommand{\cc}{\mathbb{C}}

\newcommand{\A}{\mathcal{A}}
 
\DeclareMathOperator{\gr}{gr} \DeclareMathOperator{\ad}{ad}
 
 \DeclareMathOperator{\id}{id}

\DeclareMathOperator{\Ker}{Ker} 
\DeclareMathOperator{\End}{End} 
\DeclareMathOperator{\Res}{Res} 
\DeclareMathOperator{\Spec}{Spec} 
\DeclareMathOperator{\rk}{rk} \DeclareMathOperator{\Tr}{Tr}
\renewcommand{\phi}{\varphi}

\newcommand{\ep}{\varepsilon}

\makeatletter
\def\@mult#1{\raise #1\rlap{$\cdot$}\lower #1\rlap{$\cdot$}\cdot}
\def\did{\mathrel{\@mult{3pt}}}
\makeatother
\def\openrow#1#2#3{\setbox0=\vbox{\hbox
    {\vrule height#2 width#3\kern#2\vrule height#2 width0pt}\hrule height#3}
    \hbox{\leaders\copy0\hskip#1\wd0\vrule width#3}}
\def\row#1#2#3{\vbox{\hrule height#3\openrow{#1}{#2}{#3}}}
\def\Yr#1{\row{#1}{1.5ex}{.1ex}}

\def\DY#1\endDY{\baselineskip=1ex\lineskip=0pt\lineskiplimit=0pt{\vcenter
    {\Yr#1}}}
\def\openclm#1#2#3{\setbox0=\vbox{\hrule height#3\hbox
    {\vrule width0pt\kern#2\vrule width#3 height#2}}\vtop
    {\leaders\copy0\vskip#1\ht0\hrule height#3}}
\def\clm#1#2#3{\hbox{\vrule width#3\openclm{#1}{#2}{#3}}}
\def\Yc#1{\clm{#1}{1.5ex}{.1ex}}

\def\CDY#1\endCDY{{\vcenter{\hbox{\Yc#1}}}}

\author{L.~G.~Rybnikov}
\thanks{The work was supported by CRDF grant RM1-2543, RFBR grant 05
01 00988-a and RFBR grant 05-01-02805-CNRSL-a.}
\address{Poncelet laboratory (Independent University of Moscow and CNRS) and Moscow State University,
department of Mechanics and Mathematics}
\email{leo.rybnikov@gmail.com}
\title{Argument Shift Method and Gaudin Model}

\begin{document}
\maketitle

\begin{abstract}
We construct a family of maximal commutative subalgebras in the
tensor product of $n$ copies of the universal enveloping algebra
$U(\gg)$ of a semisimple Lie algebra $\gg$. This family is
parameterized by collections $\mu; z_1,\dots,z_n$, where
$\mu\in\gg^*$, and $z_1,\dots,z_n$ are pairwise distinct complex
numbers. The construction presented here generalizes the famous
construction of the higher Gaudin hamiltonians due to Feigin,
Frenkel, and Reshetikhin. For $n=1$, our construction gives a
quantization of the family of maximal Poisson-commutative
subalgebras of $S(\gg)$ obtained by the argument shift method.
Next, we describe natural representations of commutative algebras
of our family in tensor products of finite-dimensional
$\gg$-modules as certain degenerations of the Gaudin model. In the
case of $\gg=sl_r$ we prove that our commutative subalgebras have
simple spectrum in tensor products of finite-dimensional
$\gg$-modules for generic $\mu$ and $z_i$. This implies simplicity
of spectrum in the "generic" $sl_r$ Gaudin model.
\end{abstract}

\section{Introduction}

Let $\gg$ be a semisimple complex Lie algebra, and $U(\gg)$ its
universal enveloping algebra. The algebra $U(\gg)$ bears a natural
filtration by the degree with respect to the generators. The
associated graded algebra is the symmetric algebra
$S(\gg)=\cc[\gg^*]$ by the Poincar\'e--Birkhoff--Witt theorem. The
commutator on $U(\gg)$ defines the Poisson--Lie bracket on
$S(\gg)$.

The \emph{argument shift method} gives a way to construct
subalgebras in $S(\gg)$ commutative with respect to the
Poisson--Lie bracket. The method is as follows. Let
$ZS(\gg)=S(\gg)^{\gg}$ be the center of $S(\gg)$ with respect to
the Poisson bracket, and let $\mu\in\gg^*$ be a regular semisimple
element. Then the algebra $A_{\mu}\subset S(\gg)$ generated by the
elements $\partial_{\mu}^n\Phi$, where $\Phi\in ZS(\gg)$, (or,
equivalently, generated by central elements of $S(\gg)=\cc[\gg^*]$
shifted by $t\mu$ for all $t\in\cc$) is commutative with respect
to the Poisson bracket, and has maximal possible transcendence
degree equal to $\frac{1}{2}(\dim\gg+\rk\gg)$ (see \cite{MF}).
Moreover, the subalgebras $A_{\mu}$ are maximal subalgebras in
$S(\gg)$ commutative with respect to the Poisson--Lie bracket
\cite{Tar2}. In \cite{Vin}, the subalgebras $A_{\mu}\subset
S(\gg)$ are named the {\em Mischenko--Fomenko subalgebras}.

In the present paper we lift the subalgebras $A_{\mu}\subset
S(\gg)$ to commutative subalgebras in the universal enveloping
algebra $U(\gg)$. More precisely, for any semisimple Lie algebra
$\gg$, we construct a family of commutative subalgebras
$\A_{\mu}\subset U(\gg)$ parameterized by regular semisimple
$\mu\in\gg^*$, so that $\gr\A_{\mu}=A_{\mu}$. For classical Lie
algebras $\gg$, it was done (by other methods) by Olshanski and
Nazarov (see \cite{NO, Mol}), and also by Tarasov in the case
$\gg=sl_r$ \cite{Tar1}.

The construction presented here is a modification of the famous
construction of the higher Gaudin hamiltonians (see
\cite{FFR,ER}). Gaudin model was introduced in \cite{G1} as a spin
model related to the Lie algebra $sl_2$, and generalized to the
case of an arbitrary semisimple Lie algebra in \cite{G}, 13.2.2.
The generalized Gaudin model has the following algebraic
interpretation. Let $V_{\lambda}$ be an irreducible representation
of $\gg$ with the highest weight $\lambda$. For any collection of
integral dominant weights $(\lambda)=\lambda_1,\dots,\lambda_n$,
let $V_{(\l)}=V_{\l_1}\otimes\dots\otimes V_{\l_n}$. For any
$x\in\gg$, consider the operator $x^{(i)}=1\otimes\dots\otimes
1\otimes x\otimes 1\otimes\dots\otimes 1$ ($x$ stands at the $i$th
place), acting on the space $V_{(\l)}$. Let $\{x_a\},\
a=1,\dots,\dim\gg$, be an orthonormal basis of $\gg$ with respect
to Killing form, and let $z_1,\dots,z_n$ be pairwise distinct
complex numbers. The hamiltonians of Gaudin model are the
following commuting operators acting in the space $V_{(\l)}$:
\begin{equation}\label{quadratic}
H_i=\sum\limits_{k\neq i}\sum\limits_{a=1}^{\dim\gg}
\frac{x_a^{(i)}x_a^{(k)}}{z_i-z_k}.
\end{equation}

We can regard the $H_i$ as elements of $U(\gg)^{\otimes n}$. In
\cite{FFR}, a large commutative subalgebra
$\A(z_1,\dots,z_n)\subset U(\gg)^{\otimes n}$ containing $H_i$ was
constructed. For $\gg=sl_2$, the algebra $\A(z_1,\dots,z_n)$ is
generated by $H_i$ and the central elements of $U(\gg)^{\otimes
n}$. In other cases, the algebra $\A(z_1,\dots,z_n)$ has also some
new generators known as higher Gaudin hamiltonians. The
construction of $\A(z_1,\dots,z_n)$ uses the quite nontrivial fact
\cite{FF} that the completed universal enveloping algebra of the
affine Kac--Moody algebra $\hat{\gg}$ at the critical level has a
large center $Z(\hat{\gg})$. To any collection $z_1,\dots,z_n$ of
pairwise distinct complex numbers, one can naturally assign a
homomorphism $Z(\hat{\gg})\to U(\gg)^{\otimes n}$. The image of
this homomorphism is $\A(z_1,\dots,z_n)$.

In the present paper we construct a family of homomorphisms
$Z(\hat{\gg})\to U(\gg)^{\otimes n}\otimes S(\gg)$ parameterized
by collections $z_1,\dots,z_n$ of pairwise distinct complex
numbers. For any collection $z_1,\dots,z_n$, the image of such
homomorphism is a certain commutative subalgebra
$\A(z_1,\dots,z_n,\infty)\subset U(\gg)^{\otimes n}\otimes
S(\gg)$. Taking value at any point $\mu\in\gg^*=\Spec S(\gg)$, we
obtain a commutative subalgebra $\A_{\mu}(z_1,\dots,z_n)\subset
U(\gg)^{\otimes n}$ depending on $z_1,\dots,z_n$ and
$\mu\in\gg^*$. For $n=1$, we obtain commutative subalgebras
$\A_{\mu}(z_1)=\A_{\mu}\subset U(\gg)$ which do not depend on
$z_1$. We show that $\gr\A_{\mu}=A_{\mu}$ for regular semisimple
$\mu$, i.e., the subalgebras $\A_{\mu}\subset U(\gg)$ are liftings
of Mischenko--Fomenko subalgebras (in the case of $\gg=sl_r$ this
can be deduced from Talalaev's formula for higher Gaudin
hamiltonians, cf. \cite{Tal, ChT}). For $\mu=0$, we have
$\A_0(z_1,\dots,z_n)=\A(z_1,\dots,z_n)$, i.e., the subalgebras
$\A_0(z_1,\dots,z_n)\subset U(\gg)^{\otimes n}$ are generated by
(higher) Gaudin hamiltonians. We show that the subalgebras
$\A_{\mu}(z_1,\dots,z_n)$ for generic $z_1,\dots,z_n$ and $\mu$
have maximal possible transcendence degree. These subalgebras
contain the following "non-homogeneous Gaudin hamiltonians":
$$
H_i=\sum\limits_{k\neq i}\sum\limits_{a=1}^{\dim\gg}
\frac{x_a^{(i)}x_a^{(k)}}{z_i-z_k}+\sum\limits_{a=1}^{\dim\gg}\mu(x_a)x_a^{(i)}.
$$

The main problem in Gaudin model is the problem of simultaneous
diagonalization of (higher) Gaudin hamiltonians. The bibliography
on this problem is enormous (cf. \cite{Fr1,Fr2,FFR,MV}). It
follows from the \cite{FFR} construction that all elements of
$\A(z_1,\dots,z_n)\subset U(\gg)^{\otimes n}$ are invariant with
respect to the diagonal action of $\gg$, and therefore it is
sufficient to diagonalize the algebra $\A(z_1,\dots,z_n)$ in the
subspace $V_{(\l)}^{sing}\subset V_{(\l)}$ of singular vectors
with respect to $diag_n(\gg)$ (i.e., with respect to the diagonal
action of $\gg$). The standard conjecture says that generic $z_i$
the algebra $\A(z_1,\dots,z_n)$ has simple spectrum in
$V_{(\l)}^{sing}$. This conjecture is proved in \cite{MV} for
$\gg=sl_r$ and $\l_i$ equal to $\omega_1$ or $\omega_{r-1}$ (i.e.,
for the case when every $V_{\l_i}$ is the standard representation
of $sl_r$ or its dual) and in \cite{SV} for $\gg=sl_2$ and
arbitrary $\l_i$.

It is also natural to set up a problem of diagonalization of
$\A_{\mu}(z_1,\dots,z_n)$ in the space $V_{(\l)}$. We show that
the representation of the algebra $\A_{\mu}(z_1,\dots,z_n)$ in the
space $V_{(\l)}$ is a limit of the representations of
$\A(z_1,\dots,z_{n+1})$ in $[V_{(\l)}\otimes
M^*_{z_{n+1}\mu}]^{sing}$ as $z_{n+1}\to\infty$. Here
$M^*_{z_{n+1}\mu}$ is the contragredient module of the Verma
module with highest weight $z_{n+1}\mu$, and the space
$[V_{(\l)}\otimes M^*_{z_{n+1}\mu}]^{sing}$ consists of all
singular vectors in $V_{(\l)}\otimes M^*_{z_{n+1}\mu}$ with
respect to $diag_{n+1}(\gg)$. This means that the representation
of $\A_{\mu}(z_1,\dots,z_n)$ in $V_{(\l)}$ is in some sense a
limit case of Gaudin model.

We prove the conjecture on the simplicity of the spectrum for the
representation of $\A_{\mu}(z_1,\dots,z_n)$ in the space
$V_{(\l)}$ for $\gg=sl_r$. The point of our proof is the fact that
the closure of the family $\A_{\mu}$ contains the Gelfand--Tsetlin
subalgebra (on the level of Poisson algebras, this fact was proved
by Vinberg \cite{Vin}). Hence, for $\gg=sl_r$, we conclude that
the algebra $\A_{\mu}(z_1,\dots,z_n)$ for generic $\mu$ and
$z_1,\dots,z_n$ has simple spectrum in $V_{(\l)}$ for any
$V_{(\l)}$. As a consequence, we obtain that the spectrum of the
algebra $\A_0(z_1,\dots,z_n)$ in $V_{(\l)}^{sing}$ is simple for
generic $z_i$ and $(\l)$.

The paper is organized as follows. In sections~2~and~3 we collect
some well-known facts on Mischenko--Fomenko subalgebras and the
center $Z(\hat{\gg})$ at the critical level, respectively. In
section~4 we describe the construction of the subalgebras
$\A_{\mu}$ and prove that $\gr\A_{\mu}=A_{\mu}$. In section~5 we
describe the general construction of the subalgebras
$\A_{\mu}(z_1,\dots,z_n)\subset U(\gg)^{\otimes n}$ and prove that
these subalgebras have the maximal possible transcendence degree.
In section~6 we describe the representation of
$\A_{\mu}(z_1,\dots,z_n)$ in $V_{(\l)}$ as a ``limit'' Gaudin
model. And, finally, in section~7 we prove the assertions on
simplicity of spectrum for $\gg=sl_r$.

I thank B.~L.~Feigin, E.~B.~Vinberg, V.~V.~Shuvalov,
A.~V.~Chervov, and D.~V.~Talalaev for useful discussions.

\section{Argument shift method}
Argument shift method is a particular case of the famous
Magri--Lenart construction \cite{Ma}. Let $R$ be a commutative
algebra equipped with two compatible Poisson brackets,
$\{\cdot,\cdot\}_1$ and $\{\cdot,\cdot\}_2$, (i.e., any linear
combination of $\{\cdot,\cdot\}_1$ and $\{\cdot,\cdot\}_2$ is a
Poisson bracket). Let $Z_t$ be the Poisson center of $R$ with
respect to $\{\cdot,\cdot\}_1+t\{\cdot,\cdot\}_2$. Let $A$ be the
subalgebra of $R$ generated by all $Z_t$ for generic $t$.

\begin{Fact} (cf. \cite{BB}, Proposition~4) The subalgebra $A\subset R$ is commutative with respect to any Poisson bracket
$\{\cdot,\cdot\}_1+t\{\cdot,\cdot\}_2$.
\end{Fact}

\begin{proof}
Suppose $a\in Z_{t_1},\ b\in Z_{t_2}$ with $t_1\ne t_2$. The
expression $\{a,b\}_1+t\{a,b\}_2$ is linear in $t$, and, on the
other hand, it vanishes at two distinct points, $t_1$ and $t_2$.
This means that $\{a,b\}_1+t\{a,b\}_2=0$ for all $t$.

Now suppose $a,b\in Z_{t_0}$. Since $t_0$ is generic, there exists
a continuous function $a(s)$ such that $a(t_0)=a$, and for $s$ in
a certain neighborhood of $t_0$ we have $a(s)\in Z_s$. For any $s$
in a punctured neighborhood of $t_0$ we have
$\{a(s),b\}_1+t\{a(s),b\}_2=0$, and therefore
$\{a,b\}_1+t\{a,b\}_2=0$.
\end{proof}

\begin{Corollary} Suppose that
$ZS(\gg)=S(\gg)^{\gg}$ is the center of $S(\gg)$ with respect to
the Poisson bracket, and let $\mu\in\gg^*$. Then the algebra
$A_{\mu}\subset S(\gg)$ generated by the elements
$\partial_{\mu}^n\Phi$, where $\Phi\in ZS(\gg)$, (or,
equivalently, generated by central elements of $S(\gg)=\cc[\gg^*]$
shifted by $t\mu$ for all $t\in\cc$) is commutative with respect
to the Poisson bracket.
\end{Corollary}

\begin{proof}
Take the Poisson--Lie bracket as $\{\cdot,\cdot\}_1$, and the
"frozen argument" bracket as $\{\cdot,\cdot\}_2$; this means that
for the generators we have
$$
\{x,y\}_2=\mu([x,y])\quad x,y\in\gg.
$$
Then the algebra $Z_t$ is generated by central elements of
$\cc[\gg^*]=S(\gg)$ shifted by $t\mu$.
\end{proof}

Since the Lie algebra $\gg$ is semisimple we can identify $\gg$
with $\gg^*$ and write $\mu\in\gg$.

\begin{Fact}\label{mf} \cite{MF} For regular semisimple $\mu\in\gg$ the algebra $A_{\mu}$
is a free commutative subalgebra in $S(\gg)$ with
$\frac{1}{2}(\dim\gg+\rk\gg)$ generators (this means that
$A_{\mu}$ is a commutative subalgebra of maximal possible
transcendence degree). One can take the elements
$\partial_{\mu}^n\Phi_k$, $k=1,\dots,\rk\gg$,
$n=0,1,\dots,\deg\Phi_k$, where $\Phi_k$ are basic
$\gg$-invariants in $S(\gg)$, as free generators of $A_{\mu}$.
\end{Fact}

In \cite{Sh} Shuvalov described the closure of the family of
subalgebras $A_{\mu}\subset S(\gg)$ under the condition
$\mu\in\hh^{reg}$ (i.e., for regular $\mu$ in the fixed Cartan
subalgebra). In particular, the following assertion is proved
in~\cite{Sh}.

\begin{Fact}\label{shuvalov} Suppose that $\mu(t)=\mu_0+t\mu_1+t^2\mu_2+\dots\in\hh^{reg}$
for generic $t$. Set $\z_k=\bigcap\limits_{i=0}^k\z_{\gg}(\mu_i)$
(where $\z_{\gg}(\mu_i)$ is the centralizer of $\mu_i$ in $\gg$),
$\z_{-1}=\gg$. Then we have
\begin{enumerate}
\item the subalgebra $\lim\limits_{t\to0}A_{\mu(t)}\subset S(\gg)$
is generated by all elements of $S(\z_k)^{\z_k}$ and their
derivatives (of any order) along $\mu_{k+1}$ for all $k$. \item
$\lim\limits_{t\to0}A_{\mu(t)}$ is a free commutative algebra.
\end{enumerate}
\end{Fact}

This means, in particular, that the closure of the family
$A_{\mu}$ for $\gg=sl_r$ contains the Gelfand--Tsetlin algebra
(see \cite{Vin}, 6.1--6.4). We shall discuss this case in
section~7.

The following results were obtained by Tarasov.
\begin{Fact}\label{commutant1}\cite{Tar2} The subalgebras
$A_{\mu}$ and the limit subalgebras of the type
$\lim\limits_{t\to0}A_{\mu(t)}$ are maximal commutative
subalgebras, i.e., they coincide with their Poisson centralizers
in $S(\gg)$.
\end{Fact}

The {\em symmetrization map} $\sigma:S(\gg)\to U(\gg)$ is defined
by the following property:
$$
\sigma(x^k)=x^k \quad \forall\ x\in\gg,\ k=0,1,2,\dots
$$

\begin{Fact}\label{symmetrization}\cite{Tar1,Tar3} For $\gg=sl_r$, a certain system of generators of
$A_{\mu}$ and of the limit subalgebras of the type
$\lim\limits_{t\to0}A_{\mu(t)}$ can be lifted to commuting
elements of $U(\gg)$ by the symmetrization map. This gives rise to
a unique lifting of $A_{\mu}$ to the universal enveloping algebra.
\end{Fact}

\begin{Remark} The system of generators of $A_{\mu}$ of the limit subalgebras of the type
$\lim\limits_{t\to0}A_{\mu(t)}$ to be lifted by the symmetrization
is chosen explicitly in \cite{Tar1}. It is, up to proportionality,
the system of the elements $\partial_{\mu}^n\Phi_k$,
$k=1,\dots,r-1$, $n=0,1,\dots,\deg\Phi_k$, (where $\Phi_k\in
S(sl_r)^{sl_r}$ are the coefficients of the characteristic
polynomial as functions on $sl_r$) and their limits, respectively.
We shall only use that this system of generators up to
proportionality is continuous in the parameter~$\mu$.
\end{Remark}

\section{Center at the critical level}

Let $\hat{\gg}$ be the affine Kac--Moody algebra corresponding to
$\gg$. The Lie algebra $\hat{\gg}$ is a central extension of the
formal loop algebra $\gg((t))$ by an element $K$. The commutator
relations are defined as follows:
\begin{equation}
[g_1\otimes x(t),g_2\otimes y(t)]=[g_1,g_2]\otimes
x(t)y(t)+\kappa_c(g_1,g_2)\Res_{t=0}x(t)dy(t)\cdot K,
\end{equation}
where $\kappa_c$ is the invariant scalar product on $\gg$ defined
by the formula
\begin{equation}\kappa_c(g_1,g_2)=-\frac{1}{2}\Tr_{\gg}\ad(g_1)\ad(g_2).\end{equation}
Set $\hat{\gg}_+=\gg[[t]]\subset\hat{\gg}$ and
$\hat{\gg}_-=t^{-1}\gg[t^{-1}]\subset\hat{\gg}$.

Define the completion $\tilde U(\hat{\gg})$ of $U(\hat{\gg})$ as
the inverse limit of $U(\hat{\gg})/U(\hat{\gg})(t^n\gg[[t]])$,
$n>0$. The action of $\tilde U(\hat{\gg})$ is well-defined on
$\hat{\gg}$-modules from the category $\mathcal{O}^0$ (i.e.,
$\hat{\gg}$-modules on which the Lie subalgebra $\hat{\gg}_+$ acts
locally finitely). We set $\tilde U(\hat{\gg})_c=\tilde
U(\hat{\gg})/(K-1)$. This algebra acts on $\hat{\gg}$-modules of
the {\em critical level} (i.e., $\hat{\gg}$-modules on which the
element $K$ acts as unity). The name "critical" is explained by
the fact that the representation theory at this level is most
complicated. In particular, the algebra $\tilde U(\hat{\gg})_c$
has a non-trivial center $Z(\hat{\gg})$. The following fact shows
that this center is rather large.

\begin{Fact}\label{center}\cite{FF,Fr2} \begin{enumerate}\item The natural homomorphism $Z(\hat{\gg})\to
(U(\hat{\gg})/U(\hat{\gg})(\hat{\gg}_++\cc(K-1)))^{\hat{\gg}_+}$
is surjective.

\item The Poincar\'e--Birkhoff--Witt filtration on the enveloping
algebra yields a filtration on the $\hat{\gg}_+$-module
$U(\hat{\gg})/U(\hat{\gg})(\hat{\gg}_++\cc(K-1))$. We have $\gr
(U(\hat{\gg})/U(\hat{\gg})(\hat{\gg}_++\cc(K-1)))^{\hat{\gg}_+}=(S(\hat{\gg})/S(\hat{\gg})(\hat{\gg}_++\cc
K))^{\hat{\gg}_+}$ with respect to this filtration.\end{enumerate}
\end{Fact}

Now let us give an explicit description of the algebra
$(S(\hat{\gg})/S(\hat{\gg})(\hat{\gg}_++\cc K))^{\hat{\gg}_+}$.
Since $\hat{\gg}=\hat{\gg}_+\oplus\hat{\gg}_-\oplus\cc K$ as
vector spaces, every element of
$U(\hat{\gg})/U(\hat{\gg})(\hat{\gg}_++\cc(K-1))$ (respectively,
$S(\hat{\gg})/S(\hat{\gg})(\hat{\gg}_++\cc K)$) has a unique
representative in $U(\hat{\gg}_-)$ (respectively, in
$S(\hat{\gg}_-)$). Thus we obtain the following natural embeddings
\begin{equation}
(U(\hat{\gg})/U(\hat{\gg})(\hat{\gg}_++\cc(K-1)))^{\hat{\gg}_+}\hookrightarrow
U(\hat{\gg}_-)
\end{equation}
and
\begin{equation}
(S(\hat{\gg})/S(\hat{\gg})(\hat{\gg}_++\cc
K))^{\hat{\gg}_+}\hookrightarrow S(\hat{\gg}_-).
\end{equation}
Let $\A\subset U(\hat{\gg}_-)$ and $A\subset S(\hat{\gg}_-)$ be
the images of these embeddings, respectively. Consider the
following derivations of the Lie algebra $\hat{\gg}_-$:
\begin{equation}\label{der1}
\partial_t(g\otimes t^m)=mg\otimes t^{m-1}\quad\forall g\in\gg, m=-1,-2,\dots
\end{equation}
\begin{equation}\label{der2}
t\partial_t(g\otimes t^m)=mg\otimes t^{m}\quad\forall g\in\gg,
m=-1,-2,\dots
\end{equation}
The derivations (\ref{der1}), (\ref{der2}) extend to the
derivations of the associative algebras $S(\hat{\gg}_-)$ and
$U(\hat{\gg}_-)$. The derivation (\ref{der2}) induce a grading of
these algebras.

Let $i_{-1}:S(\gg)\hookrightarrow S(\hat{\gg}_-)$ be the
embedding, which maps $g\in\gg$ to $g\otimes t^{-1}$. Let
$\Phi_k,\ k=1,\dots,\rk\gg$ be the generators of the algebra of
invariants $S(\gg)^{\gg}$.

\begin{Fact}\cite{BD,Fr2,Mu} The subalgebra $A\subset S(\hat{\gg}_-)$ is freely generated by the elements
$\partial_t^n \overline{S_k}$, $k=1,\dots,\rk\gg$,
$n=0,1,2,\dots$, where $\overline{S_k}=i_{-1}(\Phi_k)$.
\end{Fact}

It follows from the Fact~\ref{center} that the generators
$\overline{S_k}$ can be lifted to the (commuting) generators of
$\A$. This means that we have

\begin{Corollary}\label{S_k}\begin{enumerate} \item There exist the homogeneous with respect to
$t\partial_t$ elements $S_k\in\A$ such that $\gr
S_k=\overline{S_k}$. \item $\A$ is a free commutative algebra
generated by $\partial_t^n S_k$, $k=1,\dots,\rk\gg$,
$n=0,1,2,\dots$.\end{enumerate}
\end{Corollary}

In the further consideration we use only the existence of the
commutative subalgebra $\A\subset U(\hat{\gg}_-)$ and its
description from the Corollary~\ref{S_k}.

\begin{Remark} No general explicit formulas for the elements
$S_k$ are known at the moment. For the quadratic Casimir element
$\Phi_1$, the corresponding element $S_1\in\A$ is obtained from
$\overline{S_1}=i_{-1}(\Phi_1)$ by the symmetrization map. For
$\gg=sl_r$ explicit formulas for $S_k$ were obtained by Talalaev
in \cite{ChT,Tal}.
\end{Remark}

\begin{Remark} The construction of the higher Gaudin hamiltonians is as
follows. The commutative subalgebra $\A(z_1,\dots,z_n)\subset
U(\gg)^{\otimes n}$ is the image of the subalgebra $\A\subset
U(\hat{\gg}_-)$ under the homomorphism $U(\hat{\gg}_-)\to
U(\gg)^{\otimes n}$ of specialization at the points
$z_1,\dots,z_n$ (see \cite{FFR,ER}). We discuss this in Section~5.
\end{Remark}

\section{Maximal commutative subalgebras in $U(\gg)$}
For any $z\ne0$, we have an evaluation homomorphism
\begin{equation}
\phi_z:U(\hat{\gg}_-)\to U(\gg),\quad g\otimes t^m\mapsto z^mg.
\end{equation}
Furthermore, there is a homomorphism
\begin{equation}
\phi_{\infty}:U(\hat{\gg}_-)\to S(\gg),\quad g\otimes
t^{-1}\mapsto g,\ g\otimes t^m\mapsto0,\ m=-2,-3,\dots
\end{equation}

Let $\Delta:U(\hat{\gg}_-)\hookrightarrow U(\hat{\gg}_-)\otimes
U(\hat{\gg}_-)$ be the comultiplication. For any $z\ne0$, we have
the following homomorphism:
\begin{equation}
\phi_{z,\infty}=(\phi_z\otimes\phi_{\infty})\circ\Delta:U(\hat{\gg}_-)\to
U(\gg)\otimes S(\gg).
\end{equation}

More explicitly,
$$\phi_{z,\infty}(g\otimes
t^m)=z^mg\otimes 1 + \delta_{-1,m}\otimes g.$$

We set
$$
\A(z,\infty)=\phi_{z,\infty}(\A)\subset U(\gg)\otimes S(\gg)
$$

\begin{Proposition}\label{generators}
The subalgebra $\A(z,\infty)$ is generated by the coefficients of
the principal part of the Laurent series for the functions
$S_k(w)=\phi_{w-z,\infty}(S_k)$ about $z$ and by the values of
these functions at $\infty$.
\end{Proposition}
\begin{proof}
Indeed, $\A(z,\infty)$ is generated by the elements
$\phi_{z,\infty}(\partial_t^n S_k)$. These elements are Taylor
coefficients of $S_k(w)=\phi_{w-z,\infty}(S_k)$ about $w=0$. Since
$S_k(w)$ has a unique pole at $z$, the Taylor coefficients of
$S_k(w)$ about $w=0$ are linear expressions in the coefficients of
the principal part of the Laurent series for the same function
about $z$ and its value at $\infty$, and vice versa.
\end{proof}

\begin{Corollary}
The subalgebra $\A(z,\infty)\subset U(\gg)\otimes S(\gg)$ does not
depend on $z$.
\end{Corollary}
\begin{proof}
Indeed, the Laurent coefficients of the fuctions
$S_k(w)=\phi_{w-z,\infty}(S_k)$ about the point $z$ and the values
of these functions at $\infty$ do not depend on $z$.
\end{proof}

Every $\mu\in\gg^*$ defines the homomorphism of "specialization at
the point $\mu$" $S(\gg)\to\cc$. We denote this homomorphism also
by $\mu$. Consider the following family of commutative subalgebras
of $U(\gg)$, which is parameterized by $\mu\in\gg^*$:
\begin{equation}
\A_{\mu}:=(\id\otimes\mu)(\A(z,\infty))\subset U(\gg).
\end{equation}

\begin{Proposition}
All elements of the subalgebra $\A_{\mu}\subset U(\gg)$ are
$\z_{\gg}(\mu)$-invariant (where $\z_{\gg}(\mu)$ is the
centralizer of $\mu$ in $\gg$).
\end{Proposition}
\begin{proof}
Indeed, we have $\A(z,\infty)\subset[U(\gg)\otimes
S(\gg)]^{\Delta(\gg)}$, and the homomorphism $\mu$ is
$\z_{\gg}(\mu)$-equivariant. Therefore $\A_{\mu}\subset
U(\gg)^{\z_{\gg}(\mu)}$.
\end{proof}

Now let us prove that the subalgebras $\A_{\mu}\subset U(\gg)$
give a quantization of the Mischenko--Fomenko subalgebras in
$S(\gg)$ obtained by the argument shift method.

\begin{Theorem}\label{shift_arg} $\gr\A_{\mu}=A_{\mu}$ for regular semisimple $\mu\in\gg^*$.
\end{Theorem}
\begin{proof}
Let $E$ be a $\gg$-invariant derivation of $U(\hat{\gg}_-)\otimes
S(\gg)$ acting on the generators as follows:
\begin{equation}
E((g\otimes x(t))\otimes1)=1\otimes g\Res_{t=0}x(t)dt,\quad
E(1\otimes g)=0 \quad\forall\ g\in\gg.
\end{equation}
In other words, we have
$$
(g\otimes t^{-m})\otimes1\mapsto \delta_{-1,m}\otimes g,\quad
1\otimes g\mapsto0 \quad\forall\ g\in\gg.
$$

\begin{Lemma}
The subalgebra $\A(z,\infty)\subset U(\gg)\otimes S(\gg)$ is
generated by the elements
$$(\phi_z\otimes\id)(E^j(S_k\otimes1))\in U(\gg)\otimes S^j(\gg).$$
\end{Lemma}
\begin{proof}
Note that
$$(\id\otimes\phi_{\infty}\circ\Delta)(\partial_t^nS_k)=(\exp
E)(\partial_t^nS_k)\in U(\hat{\gg}_-)\otimes S(\gg).$$ Since the
elements $S_k$ are homogeneous with respect to $t\partial_t$, the
mentioned elements $(\phi_z\otimes\id)(E^j(S_k\otimes1))$ are
Laurent coefficients of the function
$S_k(w)=\phi_{w-z,\infty}(S_k)=\phi_{w-z}((\exp E)(S_k\otimes1))$
about the point $w=z$. Now it remains to take advantage of
Proposition~\ref{generators}.
\end{proof}

Now let $e$ be a $\gg$-invariant derivation of $S(\gg)\otimes
S(\gg)$ acting on the generators as follows:
\begin{equation}
e(g\otimes1)=1\otimes g,\quad e(1\otimes g)=0.
\end{equation}
Clearly, for any $f\in S(\gg)$ we have $(\id\otimes\mu)\circ
e^j(f\otimes 1)=\partial^j_{\mu}f$.

Now let us note that
$$\gr(\phi_z\otimes\id)(E^j(S_k\otimes1))=z^{(-\deg\Phi_k+j)}e^j(\Phi_k\otimes1)\in
S(\gg)\otimes S^j(\gg),$$ since $\gr S_k=i_{-1}(\Phi_k)$. Hence we
have
$$
\gr(\id\otimes\mu)\circ(\phi_z\otimes\id)(E^j(S_k\otimes1))=z^{(-\deg\Phi_k+j)}\partial_{\mu}^j(\Phi_k).
$$
Since the elements $\partial_{\mu}^j(\Phi_k)$ generate $A_{\mu}$,
we have $\gr\A_{\mu}\supset A_{\mu}$. The elements
$\partial_{\mu}^j(\Phi_k)$ are algebraically independent by
Fact~\ref{mf}, and the Lemma says that the elements
$(\id\otimes\mu)\circ(\phi_z\otimes\id)(E^j(S_k\otimes1))$
generate $\A_{\mu}$. Thus $\gr\A_{\mu}=A_{\mu}$.
\end{proof}

\section{Commutative subalgebras in $U(\gg)^{\otimes n}$}

Now let us generalize our construction. Let $U(\gg)^{\otimes n}$
be the tensor product of $n$ copies of $U(\gg)$. We denote the
subspace $1\otimes\dots\otimes 1\otimes\gg\otimes
1\otimes\dots\otimes 1\subset U(\gg)^{\otimes n}$, where $\gg$
stands at the $i$th place, by $\gg^{(i)}$. Respectively, for any
$u\in U(\gg)$ we set
\begin{equation}
u^{(i)}=1\otimes\dots\otimes 1\otimes u\otimes
1\otimes\dots\otimes 1\in U(\gg)^{\otimes n}.
\end{equation}

Let $diag_n:U(\hat{\gg}_-)\hookrightarrow U(\hat{\gg}_-)^{\otimes
n}$ be the diagonal embedding. For any collection of pairwise
distinct complex numbers $z_i, i=1,\dots,n$, we have the following
homomorphism:
\begin{equation}
\phi_{z_1,\dots,z_n,\infty}=(\phi_{z_1}\otimes\dots\otimes\phi_{z_n}
\otimes\phi_{\infty})\circ diag_{n+1}:U(\hat{\gg}_-)\to
U(\gg)^{\otimes n}\otimes S(\gg).
\end{equation}

More explicitly, we have
$$\phi_{z_1,\dots,z_n,\infty}(g\otimes
t^m)=\sum\limits_{i=1}^nz_i^mg^{(i)}\otimes 1 + \delta_{-1,m}\otimes g.$$

Set
$$
\A(z_1,\dots,z_n,\infty)=\phi_{z_1,\dots,z_n,\infty}(\A)\subset
U(\gg)^{\otimes n}\otimes S(\gg)
$$

The following assertion is proved in the same way as Proposition~1
and Corollary~2.

\begin{Proposition}\label{generators2}
\begin{enumerate}
\item The subalgebras $\A(z_1,\dots,z_n,\infty)$ are generated by
the coefficients of the principal parts of the Laurent series of
the functions
$$
S_k(w;z_1,\dots,z_n)=\phi_{w-z_1,\dots,w-z_n,\infty}(S_k)
$$
at the points $z_1,\dots,z_n$, and by their values at $\infty$.
\item The subalgebras $\A(z_1,\dots,z_n,\infty)$ are stable under
simultaneous affine transformations of the parameters $z_i\mapsto
az_i+b$. \item All the elements of $\A(z_1,\dots,z_n,\infty)$ are
invariant with respect to the diagonal action of $\gg$.
\end{enumerate}
\end{Proposition}

Let us consider the following family of commutative subalgebras in
$U(\gg)^{\otimes n}$, which is parameterized by
$z_1,\dots,z_n\in\cc$ and $\mu\in\gg^*$:
\begin{equation}
\A_{\mu}(z_1,\dots,z_n):=(\id\otimes\mu)(\A(z_1,\dots,z_n,\infty))\subset
U(\gg)^{\otimes n}.
\end{equation}

Directly from Proposition~3, we obtain

\begin{Proposition}
\begin{enumerate}
\item The subalgebras $\A_{\mu}(z_1,\dots,z_n)$ are stable under
simultaneous translations of the parameters $z_i\mapsto z_i+b$.
\item All the elements of $\A_{\mu}(z_1,\dots,z_n)$ are invariant
with respect to the diagonal action of $\z_{\gg}(\mu)$.
\end{enumerate}
\end{Proposition}

\begin{Remark} The subalgebra $\A_0(z_1,\dots,z_n)\subset U(\gg)^{\otimes n}$
can be obtained as the image of the subalgebra $\A\subset
U(\hat{\gg}_-)$ under the homomorphism
$$
\phi_{z_1,\dots,z_n}=(\phi_{z_1}\otimes\dots\otimes\phi_{z_n}
)\circ diag_n:U(\hat{\gg}_-)\to U(\gg)^{\otimes n}.
$$
These subalgebras are just the subalgebras of the higher Gaudin
hamiltonians $\A(z_1,\dots,z_n)\subset U(\gg)^{\otimes n}$
introduced in \cite{FFR} (see also \cite{ER}). The quadratic
Gaudin hamiltonians (\ref{quadratic}) are linear combinations of
the elements $\phi_{z_1,\dots,z_n}(\partial_t^n \overline{S_1})$,
$n=0,1,2,\dots$.
\end{Remark}

We shall write $\A(z_1,\dots,z_n)$ instead of
$\A_0(z_1,\dots,z_n)$.

\begin{Proposition}
The subalgebras $\A_{\mu}(z_1,\dots,z_n)$ contain the following
``non-homogeneous Gaudin hamiltonians'':\sloppy
$$
H_i=\sum\limits_{k\neq i}\sum\limits_{a=1}^{\dim\gg}
\frac{x_a^{(i)}x_a^{(k)}}{z_i-z_k}+\sum\limits_{a=1}^{\dim\gg}\mu(x_a)x_a^{(i)}.
$$
\end{Proposition}

\begin{proof}
Since the element $S_1\in\A$ is the symmetrization of
$\overline{S_1}=i_{-1}(\Phi_1)$, the element $H_i$ is the
coefficient of $\frac{1}{z-z_i}$ in the expansion of
$S_1(w;z_1,\dots,z_n)=\phi_{w-z_1,\dots,w-z_n,\infty}(S_1)$ at the
point $w=z_i$. Now it remains to apply
Proposition~\ref{generators2}.
\end{proof}

The algebra $U(\gg)^{\otimes n}\otimes S(\gg)^{\otimes m}$ has an
increasing filtration by finite-dimensional spaces,
$U(\gg)^{\otimes n}\otimes S(\gg)^{\otimes
m}=\bigcup\limits_{k=0}^{\infty} (U(\gg)^{\otimes n}\otimes
S(\gg)^{\otimes m})_{(k)}$ (by degree with respect to the
generators). We define the limit $\lim\limits_{s\to\infty}B(s)$
for any one-parameter family of subalgebras $B(s)\subset
U(\gg)^{\otimes n}\otimes S(\gg)^{\otimes m}$ as
$\bigcup\limits_{k=0}^{\infty} \lim\limits_{s\to\infty}B(s)\cap
(U(\gg)^{\otimes n}\otimes S(\gg)^{\otimes m})_{(k)}$. It is clear
that the limit of a family of {\em commutative\/} subalgebras is a
commutative subalgebra. It is also clear that passage to the limit
commutes with homomorphisms of filtered algebras (in particular,
with the projection onto any factor and with finite-dimensional
representations).

\begin{Theorem}
$\lim\limits_{s\to\infty}\A_{\mu}(sz_1,\dots,sz_n)=\A_{\mu}^{(1)}\otimes\dots\otimes\A_{\mu}^{(n)}\subset
U(\gg)^{\otimes n}$ for regular semisimple $\mu\in\gg^*$.
\end{Theorem}
\begin{proof}
\begin{Lemma}\label{predel}
$\lim\limits_{z\to\infty}\phi_z=\ep$, where
$\ep:U(\hat{\gg}_-)\to\cc\cdot1\subset~U(\gg)$ is the co-unit.
\end{Lemma}
\begin{proof}
It is sufficient to check this on the generators. We have
$$
\lim\limits_{z\to\infty}\phi_z(g\otimes
t^m)=\lim\limits_{z\to\infty}z^mg=0\quad\forall\ g\in\gg,\
m=-1,-2,\dots.
$$
\end{proof}
Now let us choose the generators of $\A(sz_1,\dots,sz_n,\infty)$
as in Proposition~\ref{generators2}. The coefficients of the
Laurent expansion of $S_k(w;sz_1,\dots,sz_n)$ at any point $sz_i$
are equal to the Laurent coefficients of
$S_k(w+sz_i;sz_1,\dots,sz_n)$ at the point $0$. On the other hand,
by Lemma~\ref{predel} we have
\begin{multline*}
\lim\limits_{s\to\infty}S_k(w+sz_i;sz_1,\dots,sz_n)=
\lim\limits_{s\to\infty}\phi_{w-s(z_1-z_i),\dots,w,\dots,w-s(z_n-z_i),\infty}(S_k)=\\
=(\ep\otimes\dots\otimes\ep\otimes\phi_w\otimes\ep\otimes\dots\otimes\ep
\otimes\phi_{\infty})\circ diag_{n+1}(S_k)=S_k^{(i)}(w;0).
\end{multline*} This means that the generators of $\A(sz_1,\dots,sz_n,\infty)$ give the generators of
$\A(z_1,\infty)^{(1)}\cdot...\cdot\A(z_n,\infty)^{(n)}$ as the
limit. Hence we conclude
$$\lim\limits_{s\to\infty}\A(sz_1,\dots,sz_n,\infty)\supset\A(z_1,\infty)^{(1)}\cdot...\cdot\A(z_n,\infty)^{(n)},$$
and therefore
$$\lim\limits_{s\to\infty}\A_{\mu}(sz_1,\dots,sz_n)\supset\A_{\mu}^{(1)}\otimes\dots\otimes\A_{\mu}^{(n)}.$$
By Fact~\ref{commutant1}, the subalgebra
$\A_{\mu}^{(1)}\otimes\dots\otimes\A_{\mu}^{(n)}\subset
U(\gg)^{\otimes n}$ coincides with its own centralizer. Thus we
have
$\lim\limits_{s\to\infty}\A_{\mu}(sz_1,\dots,sz_n)=\A_{\mu}^{(1)}\otimes\dots\otimes\A_{\mu}^{(n)}$.
\end{proof}

\begin{Corollary} For generic values of the parameters, the commutative subalgebra
$\A_{\mu}(z_1,\dots,z_n)\subset U(\gg)^{\otimes n}$  has the
maximal possible transcendence degree (which is equal to
$\frac{n}{2}(\dim\gg+\rk\gg)$).
\end{Corollary}
\begin{proof}
Indeed, for generic $\mu$ the subalgebra
$\A_{\mu}^{(1)}\otimes\dots\otimes\A_{\mu}^{(n)}\subset
U(\gg)^{\otimes n}$ has the maximal possible transcendence degree
due to the Fact~\ref{mf}. Since those subalgebras are contained in
the closure of the family $\A_{\mu}(z_1,\dots,z_n)$, the
subalgebra $\A_{\mu}(z_1,\dots,z_n)$ for generic values of the
parameters has the maximal possible transcendence degree as well.
\end{proof}

Consider the one-parameter family $U(\gg)_t$ of associative
algebras whose space of generators is $\gg$ and the defining
relations are as follows:
\begin{equation}
xy-yx=t[x,y]\quad\forall\ x,y\in\gg.
\end{equation}
For any $t\neq 0$, the map $\gg\to\gg,\ x\mapsto t^{-1}x$ induces
the associative algebra homomorphism
\begin{equation}
\psi_t:U(\gg)\tilde{\to}U(\gg)_t.
\end{equation}
For $t=0$, we have $U(\gg)_0=S(\gg)$.

Consider the commutative subalgebra $$(\id^{\otimes
n}\otimes\psi_{z^{-1}})(\A(z_1,\dots,z_n,z))\subset
U(\gg)^{\otimes n}\otimes U(\gg)_{z^{-1}}.
$$
Passing to the limit as $z\to\infty$, we obtain a certain
commutative subalgebra in $U(\gg)^{\otimes n}\otimes S(\gg)$.

\begin{Theorem}\label{predel1} $$\lim\limits_{z\to\infty}(\id^{\otimes
n}\otimes\psi_{z^{-1}})(\A(z_1,\dots,z_n,z))=\A(z_1,\dots,z_n,\infty)\subset
U(\gg)^{\otimes n}\otimes S(\gg).$$
\end{Theorem}

\begin{proof}
\begin{Lemma}\label{predel2}
$\lim\limits_{z\to\infty}\psi_{z^{-1}}\circ\phi_z=\phi_{\infty}$.
\end{Lemma}
\begin{proof}
It suffices to check this on the generators. We have
$$
\psi_{z^{-1}}\circ\phi_z(g\otimes t^m)=z\cdot z^mg\in
U(\gg)_{z^{-1}}\quad\forall\ g\in\gg,\ m=-1,-2,\dots.
$$
Hence,
$$
\lim\limits_{z\to\infty}\psi_{z^{-1}}\circ\phi_z(g\otimes
t^m)=\delta_{-1,m}g=\psi_{\infty}(g\otimes t^m)\in S(\gg).
$$
\end{proof}
Using Lemma~\ref{predel2}, we obtain
\begin{multline*}\lim\limits_{z\to\infty}(\id^{\otimes
n}\otimes\psi_{z^{-1}})(\A(z_1,\dots,z_n,z))=\lim\limits_{z\to\infty}(\phi_{z_1}\otimes\dots\otimes\phi_{z_n}
\otimes(\psi_{z^{-1}}\circ\phi_z))\circ
diag_{n+1}(\A)=\\=(\phi_{z_1}\otimes\dots\otimes\phi_{z_n}
\otimes\phi_{\infty})\circ
diag_{n+1}(\A)=\A(z_1,\dots,z_n,\infty)\subset U(\gg)^{\otimes
n}\otimes S(\gg).
\end{multline*}
\end{proof}

\section{"Limit" Gaudin model}

Let $V_{\l}$ be a finite-dimensional irreducible $\gg$-module of
the highest weight $\l$.

We consider the following $U(\gg)^{\otimes n}$-module:
\begin{equation}
V_{(\l)}:=V_{\l_1}\otimes\dots\otimes V_{\l_n}.
\end{equation}

The subalgebra $\A(z_1,\dots,z_n)\subset U(\gg)^{\otimes n}$
consists of $diag_n(\gg)$-invariant elements, and therefore acts
on the space $V_{(\l)}^{sing}\subset V_{(\l)}$ of singular vectors
with respect to $diag_n(\gg)$. This representation of
$\A(z_1,\dots,z_n)$ is known as the ($n$-point) {\em Gaudin
model}.

We will show that the representation of the subalgebra
$\A_{\mu}(z_1,\dots,z_n)\subset U(\gg)^{\otimes n}$ in the space
$V_{(\l)}$ for semisimple $\mu\in\gg^*$ is a limit case of the
$(n+1)$-point Gaudin model.

Let $M^*_{\chi}$ be the contragredient module of the Verma module
with highest weight $\chi$. This module can be constructed as
follows. Suppose that $\Delta_+$ is the set of positive roots of
$\gg$. Then we have
$M^*_{\chi}=\cc[x_{\alpha}]_{\alpha\in\Delta_+}$ (the generators
$x_{\alpha}$ have (multi-) degree $\alpha$), and the elements of
$\gg$ act by the following formulas:
\begin{enumerate}
\item The elements $e_{\alpha},\ \alpha\in\Delta_+$ of the
subalgebra $\nn_+$ act as
$$
\frac{\partial}{\partial x_{\alpha}}+\sum\limits_{\beta
> \alpha} P_{\beta}^{\alpha}\frac{\partial}{\partial x_{\beta}},
$$ where $P_{\beta}^{\alpha}$ is a certain polynomial of degree
$\beta-\alpha$. \item The elements $h\in\hh$ act as
$$
h=\chi(h)-\sum\limits_{\beta\in\Delta_+}
\beta(h)x_{\beta}\frac{\partial}{\partial x_{\beta}}.
$$
\item The generators $e_{-\alpha_i}$ (where $\alpha_i$ are the
simple roots) of the subalgebra $\nn_-$ act as
$$
e_{-\alpha_i}=\chi(h_{\alpha_i})x_{\alpha_i}+\sum\limits_{\beta\in\Delta_+}
Q_{\beta}^{\alpha_i}\frac{\partial}{\partial x_{\beta}},
$$ where $Q_{\beta}^{\alpha_i}$ is a certain polynomial of degree
$\beta+\alpha_i$.
\end{enumerate}

Consider the $U(\gg)^{\otimes n}\otimes U(\gg)$-module
$V_{(\l)}\otimes M^*_{z\mu}$. We identify the vector space
$M^*_{z\mu}$ with $\cc[x_{\alpha}]$ and re-scale the generators
setting $y_{\alpha}=z^{ht(\alpha)}x_{\alpha}$, where $ht(\alpha)$
stands for the height of a root $\alpha$. The formulas for the
action of the Lie algebra $\gg$ on $M^*_{z\mu}=\cc[y_{\alpha}]$
now look as follows:
$$
e_{\alpha}=z^{ht(\alpha)}\frac{\partial}{\partial
y_{\alpha}}+\sum\limits_{\beta
> \alpha} z^{ht(\alpha)}P_{\beta}^{\alpha}\frac{\partial}{\partial y_{\beta}},
$$
$$
h=z\mu(h)-\sum\limits_{\beta\in\Delta_+}
\beta(h)y_{\beta}\frac{\partial}{\partial y_{\beta}}.
$$
$$
e_{-\alpha_i}=\mu(h_{\alpha_i})y_{\alpha_i}+z^{-1}\sum\limits_{\beta\in\Delta_+}
Q_{\beta}^{\alpha_i}\frac{\partial}{\partial y_{\beta}},
$$

Thus we can assume that the basis of $V_{(\l)}\otimes
M^*_{z\mu}=V_{(\l)}\otimes \cc[y_{\alpha}]$ does not depend on
$z$, and the operators from $U(\gg)^{\otimes n}\otimes U(\gg)$ do
depend on $z$. The subspace of singular vectors $[V_{(\l)}\otimes
M^*_{z\mu}]^{sing}\subset V_{(\l)}\otimes \cc[y_{\alpha}]$ now
becomes depending on $z$ as well. Furthermore, the space
$V_{(\l)}\otimes M^*_{z\mu}=V_{(\l)}\otimes \cc[y_{\alpha}]$ is
graded by weights of the diagonal action of $\gg$, where the
homogeneous components do not depend on $z$ and have finite
dimensions. The subspace $[V_{(\l)}\otimes
M^*_{z\mu}]^{sing}\subset V_{(\l)}\otimes \cc[y_{\alpha}]$ is
contained in a finite sum of homogeneous components, and hence the
following limit is well-defined:
$\lim\limits_{z\to\infty}[V_{(\l)}\otimes
M^*_{z\mu}]^{sing}\subset V_{(\l)}\otimes \cc[y_{\alpha}]$.
Moreover, the limit of the image of $\A(z_1,\dots,z_n,z)$ in
$\End([V_{(\l)}\otimes M^*_{z\mu}]^{sing})$ as $z\to\infty$ is a
commutative subalgebra in
$\End(\lim\limits_{z\to\infty}[V_{(\l)}\otimes
M^*_{z\mu}]^{sing})$.

\begin{Theorem}\label{predel3} For $z\to\infty$ we have \begin{enumerate}\item the limit of $[V_{(\l)}\otimes
M^*_{z\mu}]^{sing}\subset V_{(\l)}\otimes \cc[y_{\alpha}]$ is
$V_{(\l)}\otimes1$; \item the limit of the image of
$\A(z_1,\dots,z_n,z)$ in $\End([V_{(\l)}\otimes
M^*_{z\mu}]^{sing})$ contains the image of
$\A_{\mu}(z_1,\dots,z_n)$ in
$\End(V_{(\l)}\otimes1)=\End(V_{(\l)})$.
\end{enumerate}
\end{Theorem}
\begin{proof}
Let us prove the first assertion. The subspace $[V_{(\l)}\otimes
M^*_{z\mu}]^{sing}\subset V_{(\l)}\otimes M^*_{z\mu}$ is the
intersection of kernels of the operators
$diag_{n+1}(e_{\alpha})=\sum\limits_{i=1}^{n+1}e_{\alpha}^{(i)}$,
$\alpha\in\Delta_+$. Clearly,
$$
\lim\limits_{z\to\infty}z^{-ht(\alpha)}
diag_{n+1}(e_{\alpha})=1^{\otimes
n}\otimes(\frac{\partial}{\partial y_{\alpha}}+\sum\limits_{\beta
> \alpha} P_{\beta}^{\alpha}\frac{\partial}{\partial y_{\beta}}).
$$
Hence,
$$
\lim\limits_{z\to\infty}[V_{(\l)}\otimes
M^*_{z\mu}]^{sing}\subset\bigcap\limits_{\alpha\in\Delta_+}\Ker
1^{\otimes n}\otimes(\frac{\partial}{\partial
y_{\alpha}}+\sum\limits_{\beta
> \alpha} P_{\beta}^{\alpha}\frac{\partial}{\partial y_{\beta}})
= V_{(\l)}\otimes1.
$$
Since $\dim[V_{(\l)}\otimes M^*_{z\mu}]^{sing}\ge\dim V_{(\l)}$,
we conclude that $\lim\limits_{z\to\infty}[V_{(\l)}\otimes
M^*_{z\mu}]^{sing}= V_{(\l)}\otimes1$.

Now let us prove the second assertion. The module $V_{(\l)}\otimes
M^*_{z\mu}=V_{(\l)}\otimes \cc[y_{\alpha}]$ can be regarded as
$U(\gg)^{\otimes n}\otimes U(\gg)_{z^{-1}}$-module with highest
weight $(\l_1,\dots,\l_n,\mu)$. Using the formulas for the action
of the Lie algebra $\gg$ on $\cc[y_{\alpha}]$, we see that
$$
\lim\limits_{z\to\infty}1\otimes\dots\otimes 1\otimes
e_{-\alpha_i}=\lim\limits_{z\to\infty}z^{-1}1\otimes\dots\otimes
1\otimes \psi_{z^{-1}}(e_{-\alpha_i})=0
$$
for any simple root $\alpha_i\in\Delta_+$.

Therefore, the subspace $V_{(\l)}\otimes 1\subset V_{(\l)}\otimes
\cc[y_{\alpha}]$ is stable with respect to the action of
$\lim\limits_{z\to\infty}U(\gg)^{\otimes n}\otimes
U(\gg)_{z^{-1}}=U(\gg)^{\otimes n}\otimes S(\gg)$. Moreover, the
algebra $1\otimes\dots\otimes 1\otimes S(\gg)$ acts on this space
through the character $\mu$. By Theorem~\ref{predel1}, we have
$$\lim\limits_{z\to\infty}(\id^{\otimes
n}\otimes\psi_{z^{-1}})(\A(z_1,\dots,z_n,z))=\A(z_1,\dots,z_n,\infty)\subset
U(\gg)^{\otimes n}\otimes S(\gg).$$ This means that the limit of
the image of $\A(z_1,\dots,z_n,z)$ in $\End([V_{(\l)}\otimes
M^*_{z\mu}]^{sing})$ contains the image of the algebra
$(\id\otimes\mu)(\A(z_1,\dots,z_n,\infty))=\A_{\mu}(z_1,\dots,z_n)$
in $\End(V_{(\l)}\otimes1)$.
\end{proof}

\section{The case of $sl_r$}

In this section we set $\gg=sl_r$.

\begin{Lemma}\label{gelfand-tsetlin} For $\gg=sl_r$ and $\mu(t)=E_{11}+tE_{22}+\dots+t^{n-1}E_{nn}$,
the limit subalgebra $\lim\limits_{t\to0}\A_{\mu(t)}$ is the
Gelfand--Tsetlin subalgebra in $U(sl_r)$.
\end{Lemma}
\begin{proof}
From Shuvalov's results (Fact~\ref{shuvalov}) it follows that the
associated graded algebra $\lim\limits_{t\to0}A_{\mu(t)}\subset
S(\gg)$ is the Gelfand--Tsetlin subalgebra in $S(\gg)$. Indeed, in
this case $\z_k$ is the Lie algebra $sl_{r-k-1}\oplus\cc^{k+1}$
consisting of all matrices $A\in sl_r$ satisfying
$$
A_{ij}=A_{ji}=0,\ i=1,\dots,k+1,\ j=1,\dots,r,\ i\ne j.
$$
The subalgebra of $S(sl_r)$ that is generated by $S(\z_k)^{\z_k}$
for all $k$ is the Gelfand--Tsetlin subalgebra.

For any $\mu$, the generators of $\A_{\mu}$ are the images of the
generators of $A_{\mu}$ under the symmetrization map
(Fact~\ref{symmetrization}). Therefore, the generators of
$\lim\limits_{t\to0}\A_{\mu(t)}\subset U(\gg)$ are the images of
the generators of $\lim\limits_{t\to0}A_{\mu(t)}\subset S(\gg)$
under the symmetrization map as well.

The uniqueness of the lifting (Fact~\ref{symmetrization}) implies
that $\lim\limits_{t\to0}\A_{\mu(t)}$ is the subalgebra in
$U(sl_r)$ generated by all elements of $ZU(\z_k)$ for all $k$,
i.e., it is the Gelfand--Tsetlin subalgebra in $U(sl_r)$.
\end{proof}

\begin{Theorem}\label{simple_spec} The algebra $\A_{\mu}(z_1,\dots,z_n)$ has
simple spectrum in $V_{(\l)}$ for generic values of the parameters $\mu$ and $z_1,\dots,z_n$.
\end{Theorem}

\begin{proof} \begin{enumerate} \item The Gelfand--Tsetlin
subalgebra in $U(sl_r)$ has simple spectrum in $V_{\l}$ for any
$\l$ -- it is a well-known classical result. \item Since the
Gelfand--Tsetlin subalgebra is a limit of $\A_{\mu}$, the algebra
$\A_{\mu}$ for generic $\mu$ has simple spectrum in $V_{\l}$ as
well. \item This means that for generic $\mu$ the subalgebra
$\A_{\mu}(z_1)^{(1)}\otimes\dots\otimes\A_{\mu}(z_n)^{(n)}$ has
simple spectrum in $V_{(\l)}$. Since the subalgebra
$\A_{\mu}(z_1)^{(1)}\otimes\dots\otimes\A_{\mu}(z_n)^{(n)}$
belongs to the closure of the family of subalgebras
$\A_{\mu}(z_1,\dots,z_n)$, for generic values of $\mu$ and $z_i$
the algebra $\A_{\mu}(z_1,\dots,z_n)$ has simple spectrum in
$V_{(\l)}$ as well.
\end{enumerate}
\end{proof}

\begin{Corollary} There exists a subset $W\subset\Lambda_+\times\dots\times\Lambda_+$, which is
Zariski dense in $\hh^*$ (where $\Lambda_+$ is the set of integral
dominant weights), such that for any $(\l)=(\l_1,\dots\l_n)\in W$
the Gaudin subalgebra $\A(z_1,\dots,z_n)$ has simple spectrum in
$V_{(\l)}^{sing}$ for generic values of the parameters
$z_1,\dots,z_n$.
\end{Corollary}
\begin{proof}
For fixed $\l_1,\dots\l_{n-1}$, the condition of non-simplicity of
the spectrum of $\A(z_1,\dots,z_n)$ in the space
$[V_{\l_1}\otimes\dots\otimes V_{\l_{n-1}}\otimes
M^*_{\l_n}]^{sing}$ is an algebraic condition on $\l_n\in\hh^*$
for any $z_1,\dots,z_n$. By Theorems~4~and~5 this condition is not
always satisfied. This means that the set of $\l_n\in\Lambda_+$
such that the spectrum of the algebra $\A(z_1,\dots,z_n)$ in
$[V_{\l_1}\otimes\dots\otimes V_{\l_{n-1}}\otimes
M^*_{\l_n}]^{sing}$ is simple for generic $z_1,\dots,z_n$, is
Zariski dense in $\hh^*$ for any collection $\l_1,\dots\l_{n-1}$.
Since $V_{\l_n}\subset M^*_{\l_n}$, the spectrum of the algebra
$\A(z_1,\dots,z_n)$ in the space $V_{(\l)}^{sing}$ is simple for
any of these collections $\l_1,\dots,\l_n$.
\end{proof}

\end{document}